\documentstyle[11pt]{article}
\begin{document}
\title{Generalized Shioda-Inose Structures on K3 Surfaces}
\author{Hur\c{s}it \"{O}nsiper \& Sinan Sert\"{o}z}
\date{}
\maketitle

\renewcommand{\thefootnote}{}
\footnotetext{AMS Subject Classification (1991): 14J28.}

This work concerns algebraic K3 surfaces admitting generalized
Shioda-Inose structures (Definition 1 below).
To generalize the classical Shioda-Inose structure ([S-I], [M]), one needs to
determine
finite groups with suitable actions both on K3 surfaces and on abelian
surfaces. To this end,
finite groups with symplectic actions on K3 surfaces were completely
determined in ([Mu2]) and ([X]) and in the latter article the configurations
of singularities on the quotients were also listed. On the complementary side,
Katsura's article ([K]) contains the classification of all finite groups
acting on abelian surfaces so as to yield generalized Kummer surfaces
(cf. [B] for related lattice theoretic discussion).

In this paper, using the results of ([K], [X])
we show that a K3 surface $X$ admitting a Shioda-Inose structure with
$G \neq {\bf Z}_{2}$ has $\rho(X) \ge 19$ in general and $\rho(X) = 20$
if $G$ is noncyclic.
We also show that on a singular K3 surface $X$,
all Shioda-Inose structures are induced by a unique
abelian surface.

Throughout the paper we will consider only algebraic K3 surfaces over {\bf C}.

Our notation will be as follows :

$A$ (resp. $X$) denotes an abelian (resp. an algebraic K3) surface.  \\
$A_{G}$ is the Kummer surface constructed from $A/G$ for a suitable finite
group $G$. \\
$K_{?}$ denotes the canonical class of ?. \\
$T_{?}$ = the transcendental lattice of ?. \\
$\rho(?)$ is the Picard number of ?. \\
We use the standard notation $A_{k}, D_{k}, E_{k}$ to denote the rational
singularities on surfaces. \\
$\vert G \vert$ denotes the order of the group $G$. \\

We begin with giving a precise definition of Shioda-Inose structures on K3
surfaces. For this, we first recall that
a generalized Kummer surface $A_{G}$ is a K3 surface which is
the minimal resolution of the quotient $A/G$ of an abelian surface $A$ by some
finite
group $G$ ([K], Definition 2.1).

{\bf Definition 1} : {\it A K3 surface $X$
admits a Shioda-Inose structure with group $G$ if $G $ acts on $X$
symplectically and the quotient $X/G$ is birational to a generalized Kummer
surface $A_{G}$.}

We note that generalized Kummer surfaces (in characteristic 0) arise
only if $G$ is isomorphic to one of the following 
groups ([K], Corollary 3.17):

${\bf Z}_{k}, k = 2,3,4,6$, \\
binary dihedral groups $Q_{8}, Q_{12}$ and \\
binary tetrahedral group $T_{24}$.

All of these possibilities occur ([K], Examples).

Comparing this list with the list of finite groups acting symplectically
on K3 surfaces ([X], Table 2), we see that all such $G$ appear as a group of
symplectic automorphisms of some K3 surface. Hence the question
of the existence of K3 surfaces admitting Shioda-Inose structures with group
$G$
makes sense for each of these groups $G$.

In the classical case of $G = {\bf Z}_{2}$, Morrison obtained the following
lattice theoretic characterization of K3 surfaces admitting Shioda-Inose
structure ([M], Corollary 6.4).

{\bf Theorem [M] :} {\it An algebraic K3 surface $X$ admits a Shioda-Inose
structure if and
only if $X$ satisfies one of the following conditions : \\
(i) $\rho(X) = 19$ or 20, \\
(ii) $\rho(X) = 18$ and $T_{X} = U \oplus T'$, \\
(iii) $\rho(X) = 17$ and $T_{X} = U^{2} \oplus T'$\\
where $U$ is the standard hyperbolic lattice.}\\

To contrast this case with the general situation, we include the following
elementary observation.

{\bf Lemma 2 :} {\it Given an abelian surface $A$, there exists a K3
surface $X$ with $\rho(X) = 16 + \rho(A)$ admitting a classical Shioda-Inose
structure induced by $A$}.

Proof :

Given $A$, we have $T_{A} \hookrightarrow U^{3}$ with
$signature(T_{A}) = (2, 4 - \rho(A))$. Therefore, taking
$\rho = 16 + \rho(A)$,
by the surjectivity of the period map for K3 surfaces
(cf. [M], Corollary 1.9 (ii)) we have a K3 surface $X$ with $\rho(X) = \rho$
and $T_{X}$ is isometric
to $T_{A}$. Applying ([M], Theorem 6.3) the conclusion follows. $\Box$ \\

We will see that for generalized Shioda-Inose structures one has
$\rho(X) \ge 19$. This bound on the Picard number follows from the
configuration
of the exceptional curves on $A_{G}$ for which we will need the following
result on the singularities of the quotient
$A/G$ for noncyclic $G$.

{\bf Proposition 3} : {\it If $G$ is a non-cyclic group
acting on an abelian surface $A$ to yield a generalized Kummer surface,
then the singularities of $A/G$ are given as follows :

$3A_{1} + 4D_{4}$  for $G = Q_{8}$, \\
$A_{1} + 2A_{2} + 3A_{3} + D_{5}$  for $G = Q_{12}$ and \\
$4A_{2} + 2A_{3} + A_{5}$ or $A_{1} + 4A_{2} + D_{4} + E_{6}$ for $G = T_{24}$
}.

Proof :

As is the case with analysis of this type, the proof is combinatorial
in essence and is quite standard (cf. [B], [K], [X]). We know that as we have
only quotient singularities, the
possible types of singularities are :

$A_{k}, k = 1, 2, 3, 5$ corresponding to stabilizer groups of type
${\bf Z}_{k}, k = 2, 3, 4, 6$ respectively and $D_{4}$
(resp. $D_{5}$, resp. $E_{6}$) corresponding to $Q_{8}$ (resp. $Q_{12}$,
resp. $T_{24}$). We index these types in this order with $i = 1,..,7$
and we let $n_{i}$ be the number of singular points of type $i$ on $A/G$.

Comparing the topological Euler characteristic of
$A - \{$fixed points of $G\}$
to that of \\
$A_{G} - \{$ exceptional curves $\}$, we obtain

$0 = \chi_{top}(A) = \vert G \vert(24 - \sum \chi_{i}n_{i}) + n$

where $\chi_{i}$ is the topological Euler characteristic of the configuration
corresponding to the singularity of type $i$ and $n$ is the total number of
fixed points of $G$ on $A$. Clearly we have
$n = \sum m_{i}n_{i}$ where $m_{i}$ is the index in $G$ of the stabilizer group
corresponding to $i$. Furthermore, as the lattice generated by (-2)-curves
on $A_{G}$ has rank $\leq 19$, in all cases we have

$n_{1}+2n_{2}+3n_{3}+5n_{4}+4n_{5}+5n_{6}+6n_{7} \leq 19$.

Using these restrictions together with the subgroup structure of each $G$,
the result follows.$\Box$\\

{\bf Corollary 4} : {\it If $X$ admits a Shioda-Inose structure with
$G \neq {\bf Z}_{2}$, then $\rho(X) \ge 19$ and $\rho(X) = 20$
if $G$ is noncyclic.}

Proof :

By ([I2], Corollary 1.2), we know that the Picard number of $X$ is equal to
the Picard number of the associated generalized Kummer surface $A_{G}$.
Therefore, if $G = {\bf Z}_{k}$ for $k = 3,4,6$, it follows from
([K], p. 17) that $\rho(X) \geq 19$. In case $G$ is noncyclic,
we apply Proposition 3 to see that $\rho(X) = 20$.$\Box$ \\

Next, we consider the variation of Shioda-Inose structures
with respect to the isogenies of abelian surfaces.

Given a K3-surface $X$ which admits a Shioda-Inose structure with group
$G$ and associated abelian surface $A$, we denote by $\pi_{A}$ (resp.
$\pi_{X}$) the rational covering map
$A \rightarrow A_{G}$ (resp. $X \rightarrow A_{G}$) into the corresponding
generalized Kummer surface $A_{G}$.

The following results follow by exactly the same proofs as in the
case of classical Shioda-Inose structures (cf. [S-I], [I2], [M]):

(1) $\pi_{A}^{*}(K_{A_{G}}) = K_{A}$ and $\pi_{X}^{*}(K_{A_{G}}) = K_{X}$,\\
(2) $\pi_{A}^{*} : T_{A_{G}} \rightarrow T_{A}$ (resp. $\pi_{X}^{*}$)
gives an isomorphism of lattices $T_{A_{G}} \cong T_{A}(n)$ \\
(resp. $T_{A_{G}} \cong T_{X}(n)$) where $n = \vert G \vert$,\\
(3) $T_{A}$ and $T_{X}$ are isometric.

Using these elementary observations we prove

{\bf Lemma 5} : {\it If $X$ is a singular K3 surface, then each and every
Shioda-Inose structure on $X$ is induced only by $A$.}

Proof :

We let $p_{A}, p_{A_{G}}, p_{X}$ denote the period maps of $A, A_{G}, X$
respectively.

{}From (1) above, it follows that the isometry
$\phi : T_{X} \rightarrow T_{A}$
satisfies $p_{A} \circ \phi = cp_{X}$ for some $c \in {\bf C}$.\\
If we have another abelian surface $A'$ inducing some Shioda-Inose structure
on $X$, with
corresponding isometry $\psi : T_{X} \rightarrow T_{A'}$ satisfying
$p_{A'} \circ \psi = c'p_{X}$ for some $c' \in {\bf C}$, then we get
an isometry $\phi \circ \psi^{-1} : T_{A'} \rightarrow T_{A}$.
As $\rho(A) = 4$,  $\phi \circ \psi^{-1}$ extends to an isometry
$\alpha : H^{2}(A', {\bf Z}) \rightarrow H^{2}(A, {\bf Z})$ ([S-M], Theorem 1
in Appendix) to give
$p_{A} \circ \alpha = p_{A} \circ (\phi \circ \psi^{-1}) = cc'^{-1}p_{A'}$,
and it follows that $A' \cong A$ or the dual $\hat A$ of $A$ ([S], Theorem 1).
This completes the proof because $X$ admits a classical Shioda-Inose structure
for which the associated abelian surface $A^{'}$ is self-dual
(Theorem [S-I]). $\Box$\\

Remark :

If $A_{1}$ and $A_{2}$ are two abelian surfaces which are
isogeneous, then we have
$T_{A_{1}} \otimes {\bf Q} \cong T_{A_{2}} \otimes {\bf Q}$.
Therefore two K3 surfaces $X_{1}, X_{2}$ are isogeneous in the sense of
([Mu1], Definition 1.8) if they admit Shioda-Inose structures
(not necessarily with the same group) induced from $A_{1}$, $A_{2}$
respectively ([Mu1], Remark 1.11). In case $X_{1}, X_{2}$ are singular K3
surfaces, the stronger form of isogeny follows from Lemma 5 using
([I1]); that is, we have rational maps $X_{1} \rightarrow X_{2}$,
$X_{2} \rightarrow X_{1}$ of finite degree.

\vspace{2cm}

{\bf References :}

[B] J. Bertin, Resaux de Kummer et surfaces K3, Invent. Math. 93,
267-284 (1988).

[I1] H. Inose, On defining equations of singular K3 surfaces and a notion
of isogeny, in Proc. Int. Symp. Alg. Geo., Kyoto 1977, Kinokuniya Books (1978),
495-502.

[I2] H. Inose, On certain Kummer surfaces which can be realized as nonsingular
quartic surfaces in ${\bf P^{3}}$, J. Fac. Sci. Univ. Tokyo Sect. IA
Math. 23 No 3, 545-560 (1976).

[K] T. Katsura, Generalized Kummer surfaces and their unirationality in
characteristic $p$, J. Fac. Sci. Univ. Tokyo, Sect. IA, Math. 34, 1-41 (1987).

[M] D. Morrison, On K3 surfaces with large Picard number, Inv. Math. 75,
105-121 (1984).

[Mu1] S. Mukai, On the moduli spaces of bundles on K3 surfaces I, in
Vector Bundles in Algebraic Geometry, Tata Institute of Fundamental
Research Studies 11, Oxford Uni. Press 1987, 341-413.

[Mu2] S. Mukai, Finite groups of automorphisms of K3 surfaces and the
Mathieu group, Inv. Math. 94, 183-221,(1988).

[S] T. Shioda, The period map of abelian surfaces, J. Fac. Sci. Uni. Tokyo,
25, 47-59 (1978).

[S-I] T. Shioda, H. Inose, On singular K3 surfaces, in Complex Analysis
and Algebraic Geometry, Iwanami Shoten (1977), 119-136.

[S-M] T. Shioda, N. Mitani, Singular abelian surfaces and binary quadratic
forms, in CLassification of Algebraic Varieties and Compact Complex
Manifolds, LNM No. 412, 259-287 (1974).

[X] G. Xiao, Galois covers between K3 surfaces, Ann. Inst. Fourier, Grenoble
46,
73-88 (1996).

\vspace{1cm}

\begin{tabular}{lll}

Department of Mathematics &~~~~~~~~~&Department of Mathematics\\

Middle East Technical University& ~~~~~~&Bilkent University\\

06531 Ankara, Turkey &~~~~~~&06533 Ankara, Turkey\\
&&\\
e-mail : hursit@rorqual.cc.metu.edu.tr&~~~~~~&e-mail :
sertoz@fen.bilkent.edu.tr\\

\end{tabular}

\end{document}